\documentclass[reqno]{amsart}

\usepackage[sort&compress, numbers]{natbib}
\usepackage{graphicx}
\usepackage{hyperref}
\hypersetup{
    colorlinks=true,
    linkcolor=blue,
    urlcolor=cyan,
    citecolor=blue
    }
\usepackage{bm}
\usepackage[T1]{fontenc}

\newcommand{\dar}{\mathcal{D}}
\newcommand{\lat}{\mathcal{R}}
\newcommand{\ten}{\mathbf}
\newcommand{\gt}{\bm}
\newcommand{\scalar}[1]{\langle#1\rangle}

\newcommand{\norm}[1]{\lVert#1\rVert}
\newcommand{\str}{\varepsilon}
\newcommand{\dd}{\mathrm{d}}
\DeclareMathOperator{\sgn}{sgn}
\DeclareMathOperator{\trace}{trace}
\renewcommand{\[}{\begin{equation}}
\renewcommand{\]}{\end{equation}}
\newcommand{\R}{\mathbb{R}}
\newcommand{\Z}{\mathbb{Z}}

\DeclareMathOperator{\adj}{adj}
\newcommand{\comment}[1]{#1}
\newcommand{\revise}[1]{#1}

\newtheorem{theorem}{\bf Theorem}
\newtheorem*{Theorem1}{Theorem~\ref{thm:1}}
\newtheorem{proposition}{\bf Proposition}
\newtheorem{example}{\bf Example}
\newtheorem{lemma}{\bf Lemma}
\newtheorem{definition}{\bf Definition}
\newtheorem{corollary}{\bf Corollary}
\newtheorem{remark}{\bf Remark}

\begin{document}

\title{How periodic surfaces bend}

\author{Hussein Nassar}
\address{Department of Mechanical and Aerospace Engineering, University of Missouri, Columbia, MO 65211, USA}
\thanks{Work supported by the NSF under CAREER award No. CMMI-2045881. Part of the work was completed at the d'Alembert Institute of Sorbonne Université in the summer of 2023. The author thanks the Institute members, Claire Lestringant in particular, for their hospitality. The author also thanks Basile Audoly (École Polytechnique) for insightful exchange. The author has no conflicts of interest. Data availability is not applicable as no data was generated.}
\email{nassarh@missouri.edu}

\begin{abstract}
A periodic surface is one that is invariant by a 2D lattice of translations. Deformation modes that stretch the lattice without stretching the surface are effective membrane modes. Deformation modes that bend the lattice without stretching the surface are effective bending modes. For periodic, piecewise smooth, simply connected surfaces, it is shown that the effective membrane modes are, in a sense, orthogonal to effective bending modes. This means that if a surface gains a membrane mode, it loses a bending mode, and conversely, in such a way that the total number of modes, membrane and bending combined, can never exceed 3. Various examples, inspired from curved-crease origami tessellations, illustrate the results.
\end{abstract}

\maketitle

\section{Introduction}
\comment{Slender structures in general and thin shells in particular prefer bending over stretching. Ideally, thin shells deform \emph{isometrically}, i.e., inextensionally~\cite{Landau1986}. This geometric insight has important consequences. For instance, in Saint-Venant's theory of torsion, the twisting of an \emph{open} thin-walled prismatic bar produces an axial deflection, a \emph{warping}, given by
\[
    w(s) = \alpha\int^s (x'y-xy')
\]
where $\alpha$ is the twisting rate, $(x(s),y(s))$ parametrizes the open section with a curvilinear coordinate $s$, and $\cdot'\equiv\dd/\dd s$. To find $w$, one typically solves stress balance for deflections of the form
\[
\dot{\ten x}(s,z) = (-\alpha z y(s),\alpha z x(s), w(s))
\]
where $z$ is the axial coordinate. Alternatively, it is possible to determine $w$ from purely geometric considerations by requiring that the thin-walled bar deform isometrically:}

\begin{proposition}\label{prop:warp}
    Let $\ten x: (s,z)\mapsto(x(s),y(s),z)$ describe an open thin-walled prismatic bar, i.e., a cylinder whose section $(x,y)$ is a simple curve. Then, $w$ is the unique warping such that $\dot{\ten x}$ is an infinitesimal isometric deformation of $\ten x$.
\end{proposition}
\begin{proof}
    It suffices to write the infinitesimal membrane strains of $\ten x$ produced by $\dot{\ten x}$ and to set them to~$0$.
\end{proof}

\begin{figure}
    \centering
    \includegraphics[width=0.7\linewidth]{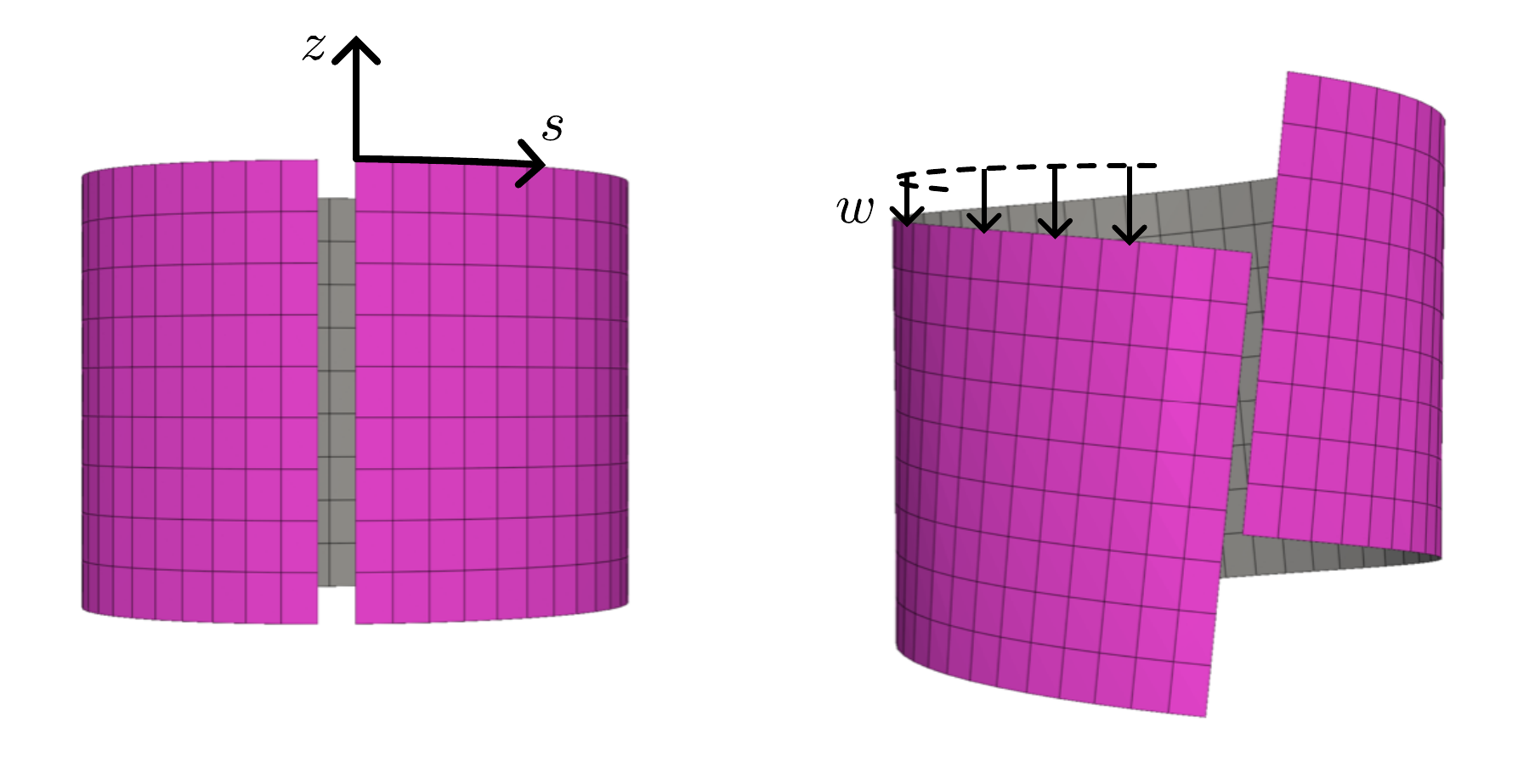}
    \caption{Warping of an open thin-walled prismatic bar before (left) and after (right) twisting.}
    \label{fig:openClose}
\end{figure}

%
\begin{remark}
    For a closed section, warping $w$ is ill-defined since it introduces a ``dislocation'' $\alpha\oint(x'y-xy') \neq 0$; see also~\cite{Landau1986b}.
\end{remark}

The above remark highlights the fact that isometric deformations, even if preferred, may not be available. In the classical mathematical literature, one finds negative results that establish the impossibility of isometric deformations for certain surfaces, often compact convex ones, e.g., Cauchy's, Dehn's, Cohn Vossen's and Pogorelov's theorems on the rigidity of convex polyhedra and surfaces~\cite{Connelly1993,Spivak1999a}. One also finds positive results that construct specific isometric deformations for specific surfaces, e.g., developable surfaces, surfaces of revolution, surfaces of translation, Cohn Vossen's surface and Connelly's flexible polyhedron~\cite{Spivak1999a,Bianchi1878,Darboux1894,Darboux1896,Nassar2023}. In the mechanics literature, isometric deformations became of interest with the birth of shell theory and, notably, the Rayleigh-Love controversy regarding the nature of dominant deformations in thin shells~\cite{Rayleigh1894,Love1906,Calladine1982}. In theory, it is now understood that the shape of the strain energy functional of a thin shell depends on whether or not its midsurface admits isometric deformations in conjunction with applied loads and boundary conditions~\cite{Harutyunyan2017,ciarlet2006}. In practice however, modern computational tools have minimized the importance of specialized geometrically-informed models (e.g., flexure shells v.s. membrane shells) and have favored more general models that, even if less efficient, can indifferently handle membrane and flexure contributions (e.g., Reissner-Mindlin theory).

Recent trends in the design, modeling and applications of compliant shell mechanisms in general and origami in particular have renewed the interest in the theory of isometric deformations~\cite{norman2009phd,Schenk2011a,Seffen2012}. In that context, much of the current literature deals with polyhedral surfaces composed of triangles or quads for which isometric deformations, sometimes referred to as ``foldings'' or ``rigid foldings'', can be constructed by solving the discrete kinematics of some planar or spherical linkages~\cite{Lang2018,Stachel2010}. A few more advance results, both positive and negative, have also been obtained for certain intrinsically flat surfaces that are creased along line and curve segments~\cite{Demaine2011,Demaine2015a,Mundilova2019}. \comment{In comparison, the aim of the present paper is to report on, and extend, a set of results regarding the availability, or impossibility, of isometric deformations for periodic surfaces, i.e., surfaces that are invariant by a 2D lattice of translations. Best known examples of such surfaces are origami and origami-like tessellations such as the Miura ori and the ``eggbox'' pattern. But other surfaces including curved-crease variants or smooth, uncreased, variants are well within the scope of the paper. Indeed, the relatively weak hypothesis of piecewise smoothness will allow to invariably handle smoothly ``corrugated'', curved-crease and polyhedral surfaces be them intrinsically flat or not. This unprecedented level of universality in the treatment is made possible by proof techniques that are free of specific constructs (e.g., spherical linkages, torsal rulings, conjugate nets) and instead use high-level arguments (e.g., symmetry, integral theorems, continuity). The other main hypothesis is that of simple connectivity: the theory excludes slits and cut-outs (e.g., kirigami).}

Much of the interest surrounding compliant shell mechanisms and origami tessellations resides in the fact that they can \emph{effectively} stretch and \emph{effectively} bend without \emph{actually} stretching~\cite{Lebee2015,Callens2017}. The main result of the proposed theory then characterizes how modes of \emph{effective} stretching of effective membrane strain $\ten E$ interact with modes of \emph{effective} bending of effective bending strain~$\gt \chi$.
\begin{Theorem1}
    Let a simply connected piecewise smooth periodic surface $\ten x$ admit an effective membrane strain $\ten E$ and an effective bending strain~$\gt \chi$, then
    \[
        E_{11}\chi_{22} - 2E_{12}\chi_{12} + E_{22}\chi_{11} = 0.
    \]
\end{Theorem1}
The theorem, quite reminiscent of a perturbative Gauss theorem~\cite{Gray2006}, \comment{establishes an orthogonality relationship between the \emph{linear spaces} of membrane modes and bending modes: the larger one space is, the smaller the other one. For instance, if the surface is free to stretch in direction~1 then it cannot bend about direction~1; if it can twist then it cannot shear, and so on. In particular, a surface can have no more than $3$ modes, bending and membrane combined (Corollary~\ref{cor:1}).}

The theorem admits another interpretation best seen when it is written in a principal basis of $\ten E$ since having $E_{12}=0$ implies
\[
    \frac{\chi_{22}}{\chi_{11}} = -\frac{E_{22}}{E_{11}}.
\]
That is: effective normal curvatures in the principal directions of effective membrane strain occur in equal and opposite proportions to the effective principal membrane strains. This is an identity between effective in-plane and out-of-plane Poisson's coefficients and, as such, has appeared and been proven for a number of periodic polyhedral surfaces with four parallelogram panels per unit cell~\cite{Schenk2013,Wei2013,Nassar2017a,Pratapa2019,Nassar2022,McInerney2022,xu2023derivation}. Theorem~\ref{thm:1} shows that in fact this identity is much more general than previously foreseen.

Two versions of Theorem~\ref{thm:1} have recently appeared in~\cite{NASSAR2024105553}, one for smooth graphs and one for a class of ``unimodal'' asymptotically isometric deformations. Here, a different version is presented for piecewise smooth surfaces in an asymptotics-free context. \comment{Beyond the proof, the main novelty resides in how Theorem~\ref{thm:1} is applied to obtain various results on the flexibility and rigidity of periodic surfaces, namely, Corollary~\ref{cor:1} and Examples~1-7.} But first, a crucial lemma of symmetry must be stated and proven.

\section{The symmetry lemma}
The purpose of this first section is to prove a property of symmetry for the differential operator of infinitesimal isometries. Basically, it is a property of symmetry of the equation $\str_{\mu\nu}=0$ albeit expressed for infinitesimal rotations rather than infinitesimal displacements. This property is not absolute and holds for a class of admissible deflections acting on periodic surfaces. Hereafter, the notions of admissibility and periodicity are respectively introduced. The lemma follows.

\subsection{Admissibility}
\begin{definition}\label{def:surf}
    A \emph{surface} is a (continuous) piecewise smooth map 
    \[
    \begin{split}
        \ten x: \R^2\supset \Omega&\to\R^3\\
        (\xi_1,\xi_2) &\mapsto \ten x(\xi_1,\xi_2)
    \end{split}
    \]
whose partial derivatives $\ten x_\alpha\equiv\partial\ten x/\partial\xi_\alpha$ are linearly independent wherever they are defined.
\end{definition}
The adopted definition is admittedly reductive. In practical terms, questions regarding self-contact and self-intersection are ignored and multi-charted surfaces are disregarded as a technical, non-essential, complication. On the plus side, a surface can be smooth or creased, where crease lines are lines across which the tangent plane experiences a jump. Also, the partials of a tensor-valued field such as $\ten x$ are denoted with a subscript as in $\ten x_1$ and $\ten x_{12}$. Otherwise, the subscript denotes a coordinate or a component as in $\xi_1$ and $\str_{12}$. \comment{Greek indices run over $\{1,2\}$.}

\begin{definition}\label{def:infiso}
    An \emph{admissible deflection} of a surface $\ten x$ is a (continuous) piecewise smooth field $\dot{\ten x}$ that is smooth wherever $\ten x$ is smooth. The \emph{infinitesimal strain} $\gt\str$ is then the $2\times 2$ matrix of coefficients
    \[  
    \str_{\mu\nu} \equiv \frac{1}{2}\left(\scalar{\dot{\ten x}_\mu,\ten x_\nu}+\scalar{\dot{\ten x}_\nu,\ten x_\mu}\right).
\]
An admissible deflection $\dot{\ten x}$ is an \emph{infinitesimal isometry} if $\str_{\mu\nu}=0$ in which case it is of the form
\[
    \dot{\ten x}_\mu = \ten w\wedge\ten x_\mu
\]
for some unique field of infinitesimal rotations $\ten w$.
\end{definition}
Thus, admissible deflections can have discontinuous derivatives at crease lines that produce further folding or unfolding. In particular, infinitesimal rotations are not expected to be continuous at crease lines. That being said, the continuity of the deflection and of the surface constrain jumps in rotations to be admissible in the following sense.

\begin{definition}\label{def:adm}
    A piecewise differentiable field $\ten w$ is an \emph{admissible} field of infinitesimal rotations of a surface $\ten x$ if $s\mapsto \ten w\wedge\dd\ten x/\dd s$ is single-valued for any $s\mapsto \gt\xi(s)\in\R^2$ that parametrizes a line of discontinuity in the tangent plane of~$\ten x$.
\end{definition}

It is now possible to fully characterize infinitesimal isometries using rotations instead of deflections. This will prove very convenient in the following.

\begin{lemma}
    On a simply connected domain, a piecewise differential field $\ten w$ is the field of infinitesimal rotations of an infinitesimal isometry $\dot{\ten x}$ of a surface $\ten x$ if and only if it is admissible and solves
\[
    \dar_{\ten x}\ten w \equiv \ten w_2\wedge\ten x_1 - \ten w_1\wedge\ten x_2 = \ten 0.
\]
\end{lemma}
\begin{proof}
    Suppose $\ten w$ is the field of infinitesimal rotations of an infinitesimal isometry $\dot{\ten x}$ of a surface~$\ten x$, then $\dot{\ten x}_\mu = \ten w\wedge\ten x_\mu$ implies $\dar_{\ten x}\ten w = \ten 0$ since $\dot{\ten x}_{\mu\nu}=\dot{\ten x}_{\nu\mu}$ and similarly for $\ten x$. The tangent $\dd\ten x/\dd s$ along a crease line $s\mapsto\gt\xi(s)$ is single-valued by continuity of $\ten x$. Similarly, $\dd\dot{\ten x}/\dd s$ is single-valued but $\dd\dot{\ten x}/\dd s = \ten w\wedge\dd \ten x/\dd s$ meaning that $\ten w$ is admissible.

    The reciprocal is a consequence of the Poincaré lemma for simply connected domains and is admitted here.
\end{proof}
The ``$\dar$'' in operator $\dar_{\ten x}$ is for Darboux who studied some of the properties of symmetry that infinitesimal isometries afford, e.g., $\dar_{\ten x}\ten w=\ten 0\implies \dar_{\ten w}\ten x=\ten 0$, and if $\dot{\ten x}$ is an infinitesimal isometry of $\ten x$ then so is $\ten x$ to $\dot{\ten x}$ \cite{Darboux1894,Darboux1896,sevennec2021les}. The main purpose of this section is to prove yet another property of symmetry of $\dar$, namely that it is a symmetric bilinear form acting on periodic admissible fields of rotations.

\subsection{Periodicity}
\begin{definition}
    Let $T_1,T_2>0$ and let $R = ]0,T_1[\times]0,T_2[$. A field $\tilde {\ten x}:\R^2\to\R^3$ is $R$-\emph{periodic}~if
    \[
        \tilde{\ten x}(\xi_1+mT_1,\xi_2+n T_2) = \tilde{\ten x}(\xi_1,\xi_2)
    \]
    for all $(\xi_1,\xi_2)\in\R^2$ and all integers $(m,n)$. A surface $\ten x$ is $R$-\emph{periodic} if it is of the form
    \[
        \ten x(\xi_1,\xi_2) = \xi_1\ten p_1 + \xi_2\ten p_2 + \tilde{\ten x}(\xi_1,\xi_2)
    \]
    where $\ten p_1$ and $\ten p_2$ are linearly independent and $\tilde{\ten x}$ is periodic.
\end{definition}
Periodicity is always understood in reference to a period $R$ which is why ``$\!R$-periodic'' is hereafter shortened to ``periodic''. One could also refer to $R$ as a ``unit cell''. But perhaps the unit cell better designates the image of $R$ or the image of $R$ projected over the plane $(\ten p_1,\ten p_2)$. In any case, here, period and unit cell are used interchangeably and what is meant, should it matter, should be clear from context. Note also that the definition differentiates between a periodic \emph{surface} and a periodic \emph{field}.

\begin{definition}\label{def:E}
    Let $\ten x$ be a periodic surface. An infinitesimal isometry $\dot{\ten x}$ is an \emph{effective membrane mode} if its field of infinitesimal rotations $\ten w$ is periodic and its \emph{effective membrane strain} $\ten E$ of components
    \[
        E_{\mu\nu} \equiv \frac{\scalar{\int \ten x_\mu, \int\dot{\ten x}_\nu} + \scalar{\int \ten x_\nu, \int\dot{\ten x}_\mu}}{2}
    \]
    is not zero, where $\int$ denotes the mean value over the period $R$, namely
    \[
        \int\cdot\, \equiv \frac{1}{\text{Area}(R)}\int_{R}\cdot\,\,\dd\xi_1\dd\xi_2.
    \]
\end{definition}
Note that field $\dot{\ten x}$ is not periodic for if it was, $\ten E$ would vanish. It is however ``morally'' periodic, i.e., periodic modulo a linear map as in
    \[
        \dot{\ten x}(\xi_1,\xi_2) = \xi_1\dot{\ten p}_1 + \xi_2\dot{\ten p}_2 + \dot{\tilde{\ten x}}(\xi_1,\xi_2)
    \]
where $\dot{\tilde{\ten x}}$ \emph{is} periodic. In that case, the action of $(\dot{\ten p}_1,\dot{\ten p}_2)$ on the unit cell $(\ten p_1,\ten p_2)$ defines the effective membrane strain, namely
\[
    E_{\mu\nu} = \frac{\scalar{\ten p_\mu,\dot{\ten p}_\nu}+\scalar{\ten p_\nu,\dot{\ten p}_\mu}}{2}.
\]
Note also that adding a constant to $\ten w$ amounts to rotating $(\dot{\ten p}_1,\dot{\ten p}_2)$ without changing $\ten E$. Thus, one could require $(\dot{\ten p}_1,\dot{\ten p}_2)$ be in the plane $(\ten p_1,\ten p_2)$. Then $\ten E$ describes a homogeneous deformation of that plane whereas $\dot{\tilde{\ten x}}$ is a periodic correction that is necessary to preserve lengths, infinitesimally speaking.

\begin{definition}\label{def:chi}
    Let $\ten x$ be a periodic surface. An infinitesimal isometry $\dot{\ten x}$ is an \emph{effective bending mode} if its field of infinitesimal rotations $\ten w$ is periodic modulo a linear map and its \emph{effective bending strain $\gt\chi$} of components
    \[
        \chi_{\mu\nu} \equiv \frac12\scalar{\int\ten w_\nu\wedge\ten \int\ten{x}_\mu+\int\ten w_\mu\wedge\ten \int\ten{x}_\nu,\ten n}
    \]
    is not zero, where $\ten n$ is the unit normal to $(\int \ten{x}_1,\int\ten{x}_2)$.
\end{definition}

It is worthwhile to justify, or rather motivate, the definition of the effective bending strain $\gt\chi$ adopted above. To do so convincingly, one must appeal to an asymptotic argument regarding the linear nature of $\ten w$ and $\ten x$ for $\gt\xi$ large enough relative to the unit cell dimensions $T_\alpha$, or conversely, for $T_\alpha$ small enough relative to $\gt\xi$.

\begin{proposition}\label{prop:chi}
Let $\dot{\ten x}$ be an effective bending mode of a periodic surface $\ten x$ and let $\ten w$ be its field of infinitesimal rotations. Then,
\[
    \epsilon^2\dot{\ten x}(\gt\xi/\epsilon)\xrightarrow[\epsilon\to 0]{}  \frac{1}{2}\xi_\mu\xi_\nu\int\ten w_\nu\wedge\int\ten x_\mu.
\]
\end{proposition}
\begin{proof}
By definition, $\ten w = \tilde{\ten w} + \xi_\alpha\ten W_\alpha$ for some periodic, piecewise differentiable and necessarily bounded field $\tilde{\ten w}$ and two constant vectors $\ten W_\alpha$. Clearly, $\ten W_\alpha=\int\ten w_\alpha$. Then, for $\epsilon\to 0$,
\[
\begin{aligned}
\epsilon^2\dot{\ten x}(\gt\xi/\epsilon)-\epsilon^2\dot{\ten x}(\ten 0) =&\, \epsilon\int_0^1\xi_\alpha\dot{\ten x}_\alpha(s\gt\xi/\epsilon)\dd s \\
=&\, \epsilon\xi_\alpha\int_0^1\ten w(s\gt\xi/\epsilon)\wedge\ten x_\alpha(s\gt\xi/\epsilon)\dd s \\
=&\, \xi_\alpha\xi_\beta\ten W_\beta\wedge\int_0^1s\ten x_\alpha(s\gt\xi/\epsilon)\dd s + \epsilon\xi_\alpha\int_0^1\tilde{\ten w}(s\gt\xi/\epsilon)\wedge\ten x_\alpha(s\gt\xi/\epsilon)\dd s \\
\to&\, \xi_\alpha\xi_\beta\ten W_\beta\wedge\int_0^1s\int\ten x_\alpha + 0\\
=&\, \frac{1}{2}\xi_\alpha\xi_\beta\ten W_\beta\wedge\int\ten x_\alpha,
\end{aligned}
\]
where the first limit is given by the Riemann-Lebesgue lemma and the second is due to boundedness.
\end{proof}

In other words, the effective bending strain $\gt\chi$ is the second fundamental form of a limit quadratic deflection $\ten W_\nu\wedge\ten p_\mu\xi_\mu\xi_\nu/2$ obtained for infinitely fine corrugations. One could obtain a similar characterization of the effective membrane strain $\ten E$ but this is not pursued here.

\subsection{Statement and proof}
It is time to state and prove the lemma of symmetry. Both lemma and proof are taken from~\cite{NASSAR2024105553} with very minor modifications and are reported here for completeness.

\begin{lemma}\label{lem:self}
    Let $\ten x$ be a periodic surface. Then,
    \[
        \int\scalar{\gt\omega,\dar_{\ten x}\ten w} = \int\scalar{\ten w,\dar_{\ten x}\gt\omega},
    \]
    for any $\gt\omega$ and $\ten w$ that are periodic and admissible.
\end{lemma}

\begin{proof}
    Let $\{R_i\}_{1\leq i\leq n}$ be a finite set of disjoint non-empty open connected sets such that $\ten x$ is smooth over $\bar{R}_i$ and such that $\cup_i\bar{R}_i=\bar R$, where $R=]0,T_1[\times]0,T_2[$ is the period of $\ten x$. Let $\partial R_{ij}=\bar{R}_i\cap(\bar{R}_j+\lat)$, where $\lat=T_1\Z\times T_2\Z$ is the periodicity lattice. Let $s\mapsto\gt\xi(s)$ parametrize one of these intersections and let the brackets $[\cdot]$ denote the jump in any quantity across the intersection. Then,
    \[\label{eq:jump}
        \begin{aligned}\relax
            [\scalar{\gt\omega,\ten w\wedge\dd\ten x/\dd s}]
            &= \scalar{[\gt\omega],\ten w\wedge\dd\ten x/\dd s} &&\text{since $\ten w$ is admissible}\\
            &= \scalar{\ten w,\dd\ten x/\dd s\wedge [\gt\omega]} &&\text{by permutation symmetry}\\
            &= \scalar{\ten w,[\dd\ten x/\dd s\wedge \gt\omega]} &&\text{by continuity of $\ten x$}\\
            &= 0 &&\text{since $\gt\omega$ is admissible}.
        \end{aligned}
    \]
    Now write
    \[
        \begin{aligned}
            &\int_R\scalar{\gt\omega,\dar_{\ten x}\ten w} \\
            =& \int_R\scalar{\gt\omega,\ten w_2\wedge\ten x_1 - \ten w_1\wedge\ten x_2} &&\text{by definition of $\dar$}\\
            =& \int_R\scalar{\gt\omega,(\ten w\wedge\ten x_1)_2 - (\ten w\wedge\ten x_2)_1} &&\text{by Schwarz theorem}\\
            =& \sum_{i}\oint_{\partial R_{i}}\scalar{\gt\omega,\ten w\wedge\dd\ten x/\dd s} - \int_R\left(\scalar{\gt\omega_2,\ten w\wedge\ten x_1}-\scalar{\gt\omega_1, \ten w\wedge\ten x_2}\right) &&\text{by the divergence theorem}\\
            =& \sum_{i}\oint_{\partial R_{i}}\scalar{\gt\omega,\ten w\wedge\dd\ten x/\dd s} + \int_R\scalar{\ten w,\dar_{\ten x}\gt\omega} &&\text{by permutation symmetry}\\
            =& \sum_{i<j}\int_{\partial R_{ij}}[\scalar{\gt\omega,\ten w\wedge\dd\ten x/\dd s}] + \int_R\scalar{\ten w,\dar_{\ten x}\gt\omega} &&\text{since $\partial R_i=\cup_j\partial R_{ij}$}\\
            =& \int_R\scalar{\ten w,\dar_{\ten x}\gt\omega} &&\text{by equation~\eqref{eq:jump}}.
        \end{aligned}
    \]
\end{proof}

Note that, by definition, periodic surfaces have a simply connected period $R$. This is a critical hypothesis without which the lemma fails in general. Indeed, the application of the divergence theorem would produce other boundary terms that do not necessarily vanish, not unless $\gt\omega$ and $\ten w$ were required to satisfy some specific boundary conditions. Mechanically speaking, the presence of holes introduces some boundary conditions whose material-dependent nature cannot be handled within the present purely geometric framework.

\section{Theorem~\ref{thm:1} and its implications}
Stating and proving the main result, i.e., Theorem~\ref{thm:1}, is now a straightforward algebraic matter. Various implications regarding the flexibility and rigidity of particular periodic surfaces follow.

\subsection{Proof of the theorem}
It is very tempting to apply the symmetry lemma to one effective membrane mode and one effective bending mode. The result follows.
\begin{theorem}\label{thm:1}
    Let $\ten x$ be a periodic surface. Then,
    \[
        E_{11}\chi_{22} - 2E_{12}\chi_{12} + E_{22}\chi_{11} = 0,
    \]
    for any effective membrane strain $\ten E$ and any effective bending strain~$\gt \chi$.
\end{theorem}
\begin{proof}
let $\ten w$ be the infinitesimal rotation of an effective bending mode of strain $\gt\chi$. Then, $\ten w = \tilde{\ten w} + \xi_\alpha\ten W_\alpha$ for some periodic, piecewise differentiable and admissible $\tilde{\ten w}$ and two constant vectors~$\ten W_\alpha$. Let $\gt\omega_o$ be a constant vector. Then,
\[
    \begin{aligned}
      0 &= \int \scalar{\gt\omega_o,\dar_{\ten x}\ten w} &&\text{since} \quad \dar_{\ten x}\ten w = 0 \\
      &= \int \scalar{\gt\omega_o,\dar_{\ten x}\tilde{\ten w}} + \int\scalar{\gt\omega,\ten W_2\wedge\ten x_1-\ten W_1\wedge\ten x_2} &&\text{by linearity}\\
      &= \int \scalar{\tilde{\ten w},\dar_{\ten x}\gt\omega_o} + \int\scalar{\gt\omega_o,\ten W_2\wedge\ten x_1-\ten W_1\wedge\ten x_2} &&\text{by symmetry of $\dar$}\\
      &= \scalar{\gt\omega_o,\ten W_2\wedge\int\ten x_1-\ten W_1\wedge\int\ten x_2} &&\text{since} \quad \dar_{\ten x}\gt\omega_o = 0.
      \end{aligned}
\]
Therefore,
\[\label{eq:symChi}
    \ten W_2\wedge\int \ten x_1 = 
    \ten W_1\wedge\int \ten x_2.
\]
Projecting over $\int\ten x_1$ and $\int\ten x_2$, it comes that $\ten W_1$ and $\ten W_2$ are both in the plane of $(\int \ten x_1,\int \ten x_2)$. Now, let $\gt \omega$ be the infinitesimal rotation of an effective membrane mode of strain $\ten E$. Then, by the same logic,
    \[\label{eq:selfChi}
    \begin{aligned}
      0 &= \int \scalar{\gt\omega,\dar_{\ten x}\ten w} &&\text{since} \quad \dar_{\ten x}\ten w = 0 \\
      &= \int \scalar{\gt\omega,\dar_{\ten x}\tilde{\ten w}} + \int\scalar{\gt\omega,\ten W_2\wedge\ten x_1-\ten W_1\wedge\ten x_2} &&\text{by linearity}\\
      &= \int \scalar{\tilde{\ten w},\dar_{\ten x}\gt\omega} + \int\scalar{\gt\omega,\ten W_2\wedge\ten x_1-\ten W_1\wedge\ten x_2} &&\text{by symmetry of $\dar$}\\
      &= \int\scalar{\gt\omega,\ten W_2\wedge\ten x_1-\ten W_1\wedge\ten x_2} &&\text{since} \quad \dar_{\ten x}\gt\omega = 0\\
      &= \scalar{\ten W_1,\int\gt\omega\wedge\ten x_2}-\scalar{\ten W_2,\int\gt\omega\wedge\ten x_1} &&\text{by permutation symmetry}.
    \end{aligned}
    \]
Finally, let $\ten p_\alpha = \int\ten x_\alpha$ and $\dot{\ten p}_\alpha = \int\gt\omega\wedge\ten x_\alpha$, and write
\[\begin{aligned}
&E_{22}\ten W_1\wedge\ten p_1-E_{12}\left(\ten W_1\wedge\ten p_2 + \ten W_2\wedge\ten p_1\right)+E_{11}\ten W_2\wedge\ten p_2 \\
&=\ten W_1\wedge(\scalar{\ten p_2,\dot{\ten p}_2}\ten p_1-\scalar{\ten p_1,\dot{\ten p}_2}\ten p_2) + \ten W_2\wedge(\scalar{\ten p_1,\dot{\ten p}_1}\ten p_2-\scalar{\ten p_2,\dot{\ten p}_1}\ten p_1)\\
&= \ten W_1\wedge(\dot{\ten p}_2\wedge(\ten p_1\wedge\ten p_2))-\ten W_2\wedge(\dot{\ten p}_1\wedge(\ten p_1\wedge\ten p_2))\\
&=-\scalar{\ten W_1,\dot{\ten p}_2}\ten p_1\wedge\ten p_2+\scalar{\ten W_2,\dot{\ten p}_1}\ten p_1\wedge\ten p_2\\
&= \ten 0,
\end{aligned}
\]
where the definition of $\ten E$, the symmetry~\eqref{eq:symChi}, the formula of the triple cross product to factor then to expand, the orthogonality $\ten W_\alpha\perp\ten p_1\wedge\ten p_2$ and equation~\eqref{eq:selfChi} have been used, respectively. The component parallel to $\ten p_1\wedge\ten p_2$ is the desired identity.
\end{proof}

\begin{corollary}\label{cor:1}
    Let $\ten x$ be a periodic surface. Let $\{\ten E\}$ and $\{\gt\chi\}$ be the linear spaces of effective membrane and bending strains. Then,
    \[
    \dim\{\ten E\} + \dim\{\gt \chi\} \leq 3.
    \]
\end{corollary}
\begin{proof}
    By Theorem~\ref{thm:1} and the rank-nullity theorem.
\end{proof}

\subsection{Examples}

\begin{figure}
    \centering
    \includegraphics[width=\linewidth]{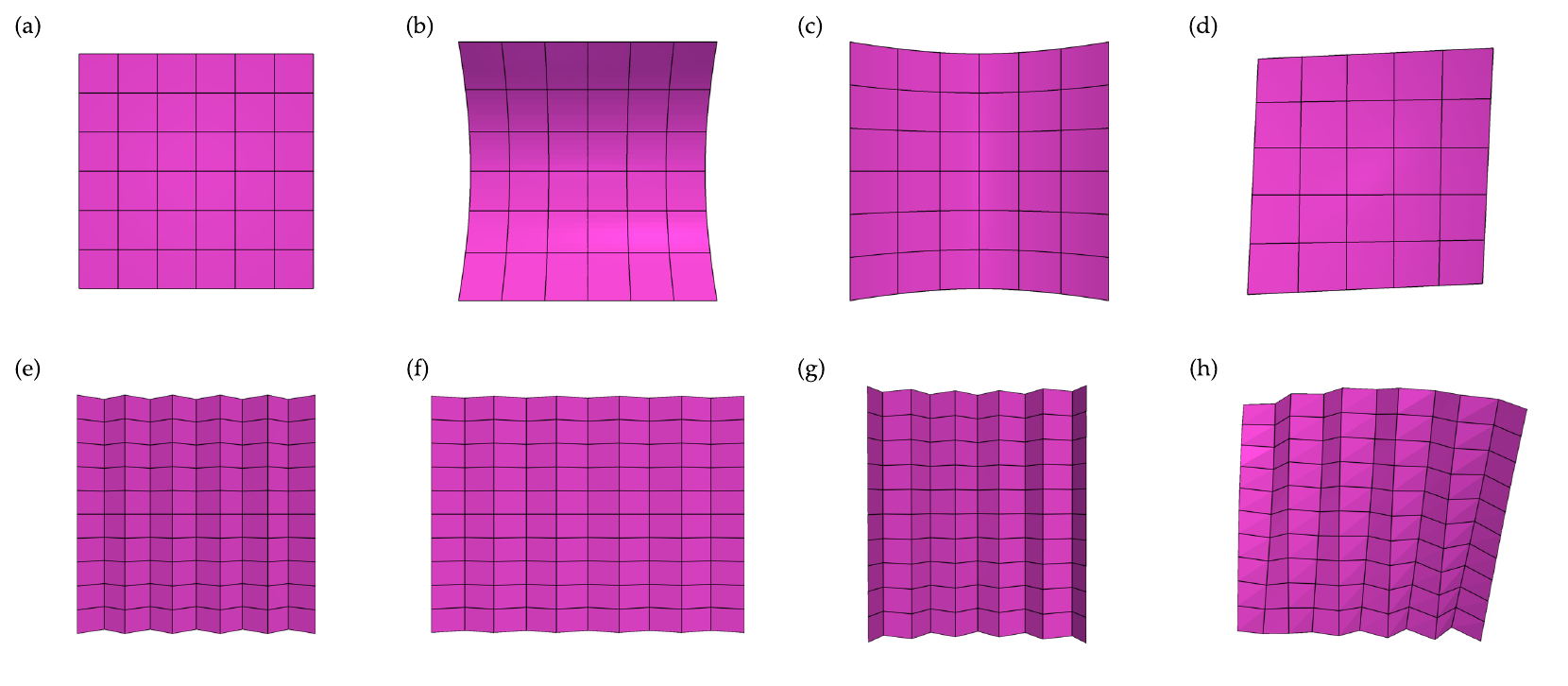}
    \caption{Orthogonality illustrated: (a) a plane and (b-d) its effective bending modes; (e) a simple corrugation, (f) its effective membrane mode and (g, h) its effective bending modes. By corrugating the plane, mode (b) is lost but mode (f) is gained.}
    \label{fig:PlanevsSimple}
\end{figure}

\begin{example}\label{ex:plane}
    The plane $\ten x:(\xi_1,\xi_2)\mapsto (\xi_1,\xi_2,0)$ is a periodic surface. For any symmetric matrix $\gt\chi$, the deflection $\dot{\ten x}:(\xi_1,\xi_2)\mapsto (0,0,\chi_{\mu\nu}\xi_\mu\xi_\nu/2)$ is an effective bending mode of effective bending strain $\gt\chi$ \revise{by Definitions~\ref{def:infiso} and~\ref{def:chi}.} Now let $\ten E$ be an effective membrane strain, then for any $\gt\chi$,
    \[
        E_{11}\chi_{22} - 2E_{12}\chi_{12} + E_{22}\chi_{11} = 0
    \]
    Hence, $\ten E = \ten 0$. In other words, the plane admits no effective membrane strains. Theorem~\ref{thm:1} appears to say that: since the plane is so flexible out of the plane, it must be completely stiff in the plane. See Figure~\ref{fig:PlanevsSimple}(a-d).
\end{example}

\begin{example}\label{ex:simpleCor}
    Let $f$ be a (continuous) piecewise smooth non-constant periodic function and let $\ten x:(\xi_1,\xi_2)\mapsto (\xi_1,\xi_2,f(\xi_1))$ be a ``simply corrugated'' periodic surface. Then, the deflection
    \[
        \dot{\ten x}:(\xi_1,\xi_2)\mapsto \left(\int^{\xi_1}f'^2,0,-f(\xi_1)\right)
    \]
    is an effective membrane mode of effective membrane strain $\ten E$ with components
    \[
        E_{11} = \int f'^2 \neq 0, \quad E_{12}=E_{22}=0
    \]
    \revise{by Definitions~\ref{def:infiso} and~\ref{def:E}.} Then, by Theorem~\ref{thm:1}, any effective bending strain $\gt\chi$ has $\chi_{22}=0$. The theorem thus maintains a trade-off between flexibility and rigidity in- and out-of-plane. Compared to the plane (Example~\ref{ex:plane}), the corrugation $f$ grants the periodic surface an effective membrane mode but takes away an effective bending mode; see Figure~\ref{fig:PlanevsSimple}(e-h). This is but a re-interpretation of Gauss theorem albeit using global constructs rather than a local one, i.e., effective modes v.s. Gaussian curvature.
\end{example}

\begin{figure}
    \centering
    \includegraphics[width=\linewidth]{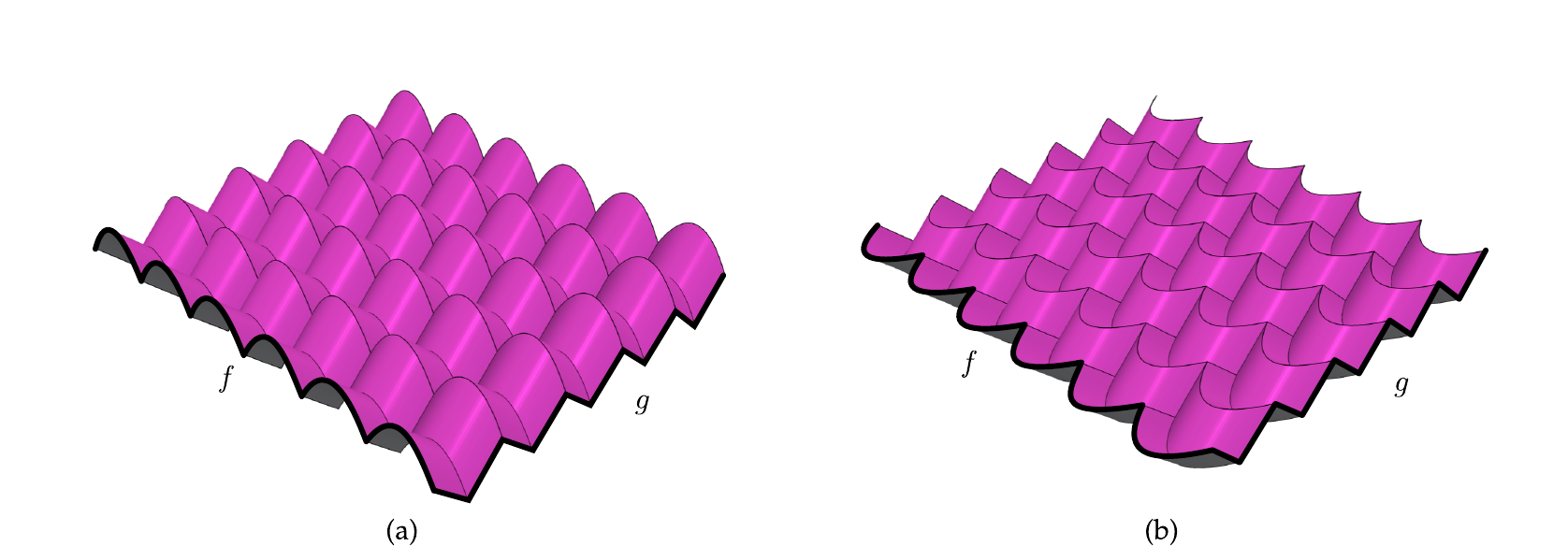}
    \caption{\revise{Two surfaces of translation: in both cases $f$ is piecewise quadratic and $g$ is piecewise linear but in case~(a), the profiles are both perpendicular to the plane of periodicity (Example~\ref{ex:eggbox} and Figure~\ref{fig:eggsvsMiura}(a-d)). By contrast, in case~(b), one profile belongs to the plane of periodicity (Example~\ref{ex:miura} and Figure~\ref{fig:eggsvsMiura}(e-h))}.}
    \label{fig:surfOfTranslation}
\end{figure}

\revise{The surfaces exemplified next are \emph{surfaces of translation}: they are obtained by translating one curve of profile $f$ along another curve of profile $g$. The construction is illustrated in Figure~\ref{fig:surfOfTranslation}.}

\begin{figure}
    \centering
    \includegraphics[width=\linewidth]{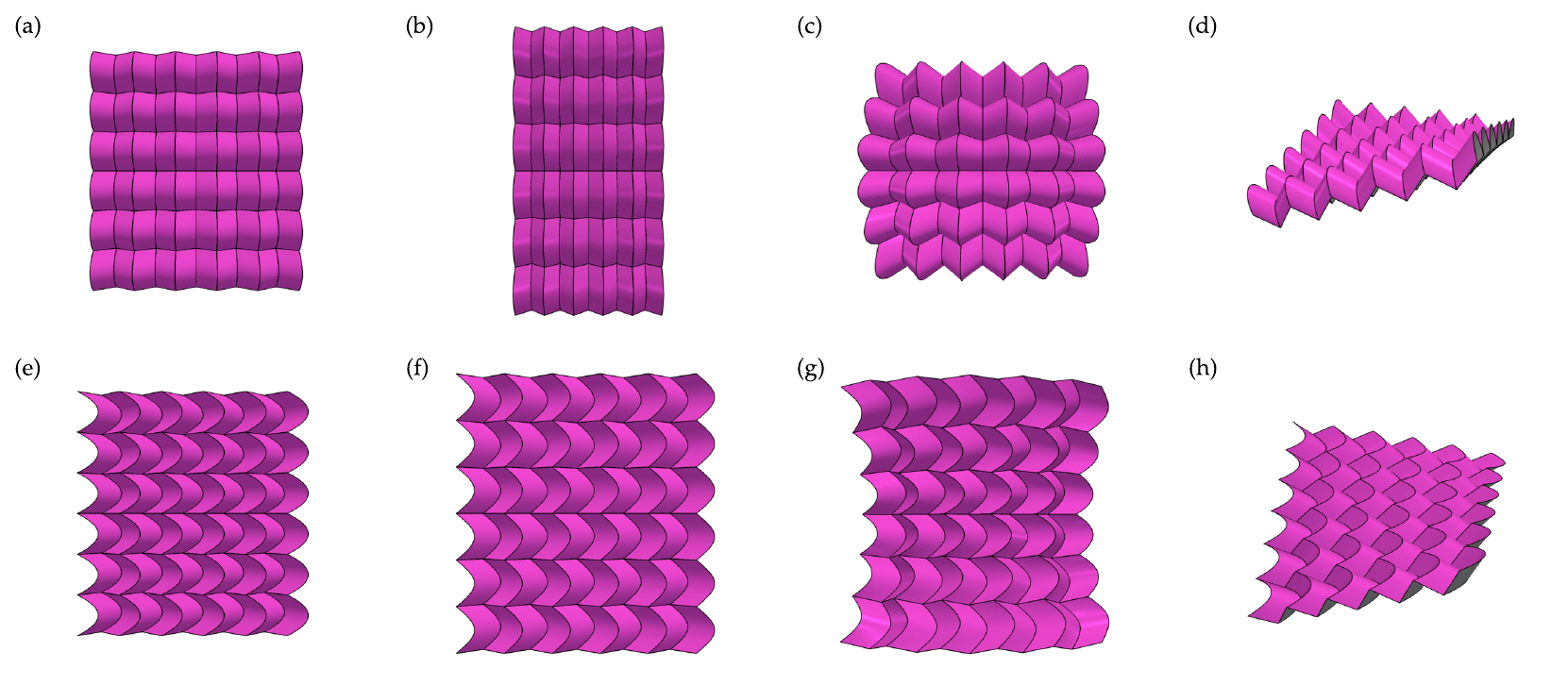}
    \caption{Orthogonality illustrated: (a) a surface from Example~\ref{ex:eggbox} where $f$ is piecewise quadratic and $g$ is piecewise linear; (b) its effective membrane mode; (c, d) its effective bending modes, (e) a surface from Example~\ref{ex:miura} with the same $f$ and $g$; (f) its effective membrane mode; (g, h) its effective bending modes. However modes (b,f) couple extensions, modes (c,g) couple curvatures in the opposite way. Modes (d, h) preclude effective shear membrane modes. Surfaces constructed by triangulation; code available online~\cite{zigzagSums}.}
    \label{fig:eggsvsMiura}
\end{figure}

\begin{example}\label{ex:eggbox}
    Let $f$ and $g$ be two (continuous) piecewise smooth non-constant periodic functions and let $\ten x:(\xi_1,\xi_2)\mapsto(\xi_1,\xi_2,f(\xi_1)+g(\xi_2))$ be a ``doubly corrugated'' periodic surface. Then, \revise{by direct verification of Definitions~\ref{def:infiso} and~\ref{def:E},} the deflection
    \[
        \dot{\ten x}:(\xi_1,\xi_2)\mapsto \left(\int^{\xi_1}f'^2,-\int^{\xi_2}g'^2,g(\xi_2)-f(\xi_1)\right)
    \]
    is an effective membrane mode of effective membrane strain
    \[
        [\ten E] = \begin{bmatrix}
            \int f'^2 & 0 \\ 0 & -\int g'^2
        \end{bmatrix}.
    \]
    Then, by Theorem~\ref{thm:1}, any effective bending $\gt\chi$ satisfies
    \[
        \frac{\chi_{22}}{\chi_{11}} = \frac{\int g'^2}{\int f'^2},
    \]
    should the ratio be defined. Thus, the double corrugation couples extension and contraction in directions~$(1,0)$ and~$(0,1)$ in the effective membrane mode and, necessarily then, couples the bending in the same directions and in the same proportions but in the opposite way.

    \revise{The ``eggbox'' pattern is a particular case where $f$ and $g$ are both piecewise linear (e.g., $f'=g'=\sgn(\cos)$). A hybrid curved-crease straight-crease variant is obtained by letting $f$ be piecewise quadratic and $g$ be piecewise linear as shown earlier on Figure~\ref{fig:surfOfTranslation}(a). The corresponding modes of deformation are shown on Figure~\ref{fig:eggsvsMiura}(a-d). As expected, the longitudinal and lateral effective membrane strains are of opposite signs, i.e., the surface stretches laterally when contracted longitudinally (panel~b). Accordingly, the effective normal curvatures are of the same sign, i.e., the surface bends into a dome (panel c).}
\end{example}

\begin{example}\label{ex:miura}
    Let $f$ and $g$ be two (continuous) piecewise smooth non-constant periodic functions and let $\ten x:(\xi_1,\xi_2)\mapsto(\xi_1,\xi_2+f(\xi_1),g(\xi_2))$. Suppose $\{g'=0\}$ is (essentially) empty. Then, \revise{by direct verification of Definitions~\ref{def:infiso} and~\ref{def:E},} the deflection
    \[
        \dot{\ten x}:(\xi_1,\xi_2)\mapsto \left(\int^{\xi_1}f'^2,-f(\xi_1)+\xi_2,-\int^{\xi_2}\frac{1}{g'}\right)
    \]
    is an effective membrane mode of effective membrane strain
    \[
        [\ten E] = \begin{bmatrix}
            \int f'^2 & 0 \\ 0 & 1
        \end{bmatrix}.
    \]
    Thus, by Theorem~\ref{thm:1}, any effective bending $\gt\chi$ satisfies
    \[
        \frac{\chi_{22}}{\chi_{11}} = -\frac{1}{\int f'^2},
    \]
    should the ratio be defined. Therefore, if $\gt\chi$ is an effective bending strain then $\det\gt\chi < 0$. In other words, for any $f$ and $g$ as stated, $\ten x$ bends ``anti-clastically'' into a saddle. See Figure~\ref{fig:eggsvsMiura}(e-h).

    \revise{The Miura ori is a particular case where $f$ and $g$ are both piecewise linear (e.g., $f'=g'=\sgn(\cos)$). Here too, a hybrid curved-crease straight-crease variant is obtained by letting $f$ be piecewise quadratic and $g$ be piecewise linear as shown earlier on Figure~\ref{fig:surfOfTranslation}(b). The corresponding modes of deformation are shown on Figure~\ref{fig:eggsvsMiura}(e-h). Indeed, the longitudinal and lateral effective membrane strains are of the same sign, i.e., the surface stretches laterally when stretched longitudinally (panel~f). Accordingly, the effective normal curvatures are of opposite signs, i.e., the surface bends into a saddle (panel g).}
\end{example}

\begin{example}\label{ex:noShear}
    The surfaces exemplified so far all admit an effective bending mode that is a pure twisting, i.e., with $\chi_{11}=\chi_{22}=0$ and $\chi_{12}\neq0$; \comment{see Figure~\ref{fig:PlanevsSimple}(d,h) and Figure~\ref{fig:eggsvsMiura}(d,h)}. This is because all of them are surfaces of translation. Here is the general case.

    \begin{proposition}
        Let $\ten x:(\xi_1,\xi_2)\mapsto\gt\alpha(\xi_1)+\gt\beta(\xi_2)$ be a periodic surface. Then, $\ten w:(\xi_1,\xi_2)\mapsto\gt\alpha(\xi_1)-\gt\beta(\xi_2)$ is the infinitesimal rotation of an effective bending mode of strain
        \[
            [\gt \chi] = \begin{bmatrix}
            0 & \norm{\int\gt\alpha'\wedge\int\gt\beta'}\\
            \norm{\int\gt\alpha'\wedge\int\gt\beta'} & 0
            \end{bmatrix}.
        \]
    \end{proposition}
    \begin{proof}
        By direct verification of Definition~\ref{def:chi}.
    \end{proof}

    Then, by Theorem~\ref{thm:1}, for any periodic surface of translation, if $\ten E$ is an effective membrane strain, then $E_{12}=0$. \revise{In other words, since these surfaces can twist, they cannot shear (relative to the same axes).}
\end{example}

At this stage, it is worthwhile to recall that Definition~\ref{def:surf} identifies surfaces and their parametrizations for convenience. That being said, the results of the theory, and Theorem~\ref{thm:1} in particular, remain meaningful if stated for a surface $\ten x(\R ^2)$ rather than for a parametrization $\ten x$. Indeed, it is possible to define the effective membrane and bending strains in a parametrization-independent fashion as one would in continuum mechanics for instance. It is equally possible to state Theorem~\ref{thm:1} using index-free notation, e.g., for any effective membrane and bending strains $\ten E$ and $\gt\chi$, one has
\[
    \lim_{t\to 0}\frac{\dd}{\dd t}\det(\ten E+t\gt\chi)=0,
\]
or
\[
    \trace(\adj(\ten E)\gt\chi)=0,
\]
where $\adj(\ten E)$ is the adjugate of $\ten E$.

\begin{example}
    Going back to Example~\ref{ex:eggbox}, reparametrize the considered surfaces using $\ten x: (\xi_1,\xi_2)\mapsto(\xi_1,\xi_2+\gamma\xi_1,f(\xi_1)+g(\xi_2+\gamma\xi_1)$ where $\gamma=T_2/T_1$ is ratio of the period of $g$ to that of $f$. Then, the deflection
    \[
        \dot{\ten x}:(\xi_1,\xi_2)\mapsto \left(\int^{\xi_1}f'^2,-\int^{\xi_2+\gamma\xi_1}g'^2,g(\xi_2+\gamma\xi_1)-f(\xi_1)\right)
    \]
    is an effective membrane mode of effective membrane strain
    \[
        [\ten E] = \begin{bmatrix}
            \int f'^2-\gamma^2\int g'^2 & -\gamma \int g'^2 \\ -\gamma\int g'^2 & -\int g'^2
        \end{bmatrix} = \begin{bmatrix}
            1 & \gamma \\ 0 & 1
        \end{bmatrix}\begin{bmatrix}
            \int f'^2 & 0 \\ 0 & -\int g'^2
        \end{bmatrix}\begin{bmatrix}
            1 & 0 \\ \gamma & 1
        \end{bmatrix},
    \]
    as one would expect by transforming the components from Example~\ref{ex:eggbox}. Similarly, by Theorem~\ref{thm:1}, any effective bending $\gt\chi$ satisfies
    \[
        \left(\int f'^2 -\gamma^2\int g'^2\right)\chi_{22}-2\gamma\int g'^2 \chi_{12} -\int g'^2\chi_{11} = 0,
    \]
    which can be rearranged into
    \[
        \int f'^2\chi_{22} - \int g'^2(\gamma^2\chi_{22}+2\gamma\chi_{12}+\chi_{11}) = 0,
    \]
    again, as one would expect from Example~\ref{ex:eggbox} by transforming the components of $\gt\chi$ like a bilinear form.
\end{example}

\begin{example}
    Let $\ten x$ be a periodic surface that admits an effective membrane strain $\ten E$. The membrane strain being symmetric, there exists an orthonormal basis in which its matrix is diagonal. Then, in that basis, by Theorem~\ref{thm:1},
    \[
        E_{11}\chi_{22} + E_{22}\chi_{11} = 0,
    \]
    for any effective bending strain $\gt\chi$. Re-arrange, if possible, into
    \[\label{eq:identity}
        \frac{E_{22}}{E_{11}} = -\frac{\chi_{22}}{\chi_{11}}
    \]
    to deduce that: for any (piecewise smooth, simply connected) periodic surface, the ratio of effective principal membrane strains is equal and opposite to the ratio of effective normal curvatures in the principal directions of effective membrane strain.
\end{example}

\subsection{Further discussion}
Identity~\eqref{eq:identity} has been proven and verified numerically in a number of particular cases~\cite{Schenk2013,Wei2013,Nassar2017a,Pratapa2019,Nassar2022,McInerney2022,xu2023derivation,Karami2024}. There is some confusion however regarding interpretation and that warrants further clarification. Suppose that the effective membrane strain has $E_{12}=0$. In that case, Identity~\eqref{eq:identity} holds albeit in a basis that is not necessarily orthonormal. But then again, identity~\eqref{eq:identity} also holds in an orthonormal basis aligned with the principal directions of $\ten E$. Both things can be true but perhaps the term ``Poisson's coefficient'' is better reserved for the value that $-E_{22}/E_{11}$ takes in an orthonormal basis.

It is seen that Theorem~\ref{thm:1} is mainly used in two ways. Either it leverages the existence of some effective membrane modes to eliminate certain effective bending modes or it leverages the existence of some effective bending modes to eliminate certain effective membrane modes. It does so as if to preserve a measure of flexibility. Now, Corollary~\ref{cor:1} ensures that
\[
    \dim\{\ten E\} + \dim\{\gt\chi\} \leq 3,
\]
but in all of the above examples, the equality holds. One could probably conjure examples where the strict inequality holds (e.g., a Miura ori interspersed with thin flat strips) but such ``intentional counter-examples'' are not pursued here. Harder to produce are examples where the number of effective membrane modes, i.e., $\dim\{\ten E\}$, exceeds that of effective bending modes, i.e., $\dim\{\gt\chi\}$. More importantly, it is worthwhile to recall that the topology of the surface, its simple connectedness in particular, is a main ingredient of the theory. Should the surface have holes or handles, integrability becomes more demanding and the extra integrability requirements provide extra rigidity and bring down the number of effective modes as was the case in the introductory example of Proposition~\ref{prop:warp}. Conversely, if the surface is not path-connected (e.g., a lattice of spheres), then the theory fails and $\dim\{\ten E\} + \dim\{\gt\chi\}$ can be trivially as high as $6$.

There is in fact one other way in which Theorem~\ref{thm:1} can be useful and that is in the spirit of Proposition~\ref{prop:chi}. Consider for instance the case of a periodic surface that admits a unique effective membrane strain $\ten E$. Then, one can claim that the infinitesimal isometries of $\ten x^\epsilon:\gt\xi\mapsto \epsilon\ten x(\gt\xi/\epsilon)$, in the limit $\epsilon\to 0$, produce perturbations $\dot{\ten X}$ to the metric of the plane $\ten X \equiv \lim_{\epsilon\to0}\ten x^\epsilon$ such that 
\[
        \frac{\scalar{\dot{\ten X}_\mu,\ten X_\nu}+\scalar{\dot{\ten X}_\nu,\ten X_\mu}}{2} = a E_{\mu\nu},
\]
where $a:\gt\xi\mapsto a(\gt\xi)$ is a scalar field that controls the amplitude of the effective membrane mode. Then, by Theorem~\ref{thm:1}, the correction $\gt\chi$ to the second fundamental form of $\ten X$ is to be found within the linear space defined by
\[
E_{11}\chi_{22} - 2E_{12}\chi_{12} + E_{22}\chi_{11} = 0.
\]
Such an asymptotic description of the isometries of a periodic surface has been successful in predicting the folded shapes of several origami tessellations, see, e.g., \cite{Nassar2018b,marazzato2024}.

\section{Conclusion}
How do periodic surfaces bend then? The proposed theory does not provide a direct answer. Instead, it characterizes how effective membrane modes and effective bending modes interact and shape each other through an orthogonality relationship:
\[
\forall\ten E,\gt\chi, \quad \trace(\adj(\ten E)\gt\chi)=0.
\]
Thus, by gaining an effective membrane mode, a periodic surface loses an effective bending mode so that the number of independent modes, membrane and bending combined, can never exceed~3.

The theory makes a certain number of assumptions with the main one being that of simple connectivity. The relationship between topology and rigidity is thematic of many structural problems as illustrated in the introductory example or in Saint-Venant's theory of torsion more generally and has been leveraged in the context of origami structures in particular; see, e.g., \cite{Filipov2015}. The present theory hopefully provides a new appreciation of how topology can contribute to geometric rigidity, as well as to elastic stiffness.

The proposed theory is purely geometric and its relevance to the behavior of elastic shells is limited to thin shells. On that front, preliminary finite element simulations suggest that as the thickness of a periodic shell is reduced, the predictions of the theory become more accurate, see, e.g.,~\cite{NASSAR2024105553}.

The techniques used in the proofs are believed to be new to the field of origami and compliant shell mechanisms and rely on some integral identities that are indifferent to smoothness or developability assumptions. They do rely on something however and the crucial symmetry lemma relies on periodicity, or at least ``closure'' in the sense that the application of the divergence theorem does not produce ``loose'' boundaries. This is quite reminiscent, albeit under stronger smoothness hypothesis, of the proof of rigidity of smooth convex compact surfaces that uses a certain integral formula of Blaschke, see~\cite{Spivak1999a}. In reality, the symmetry lemma is slightly stronger than stated in Lemma~\ref{lem:self} for it also applies under boundary conditions of periodicity modulo a rotation. These are, for instance, the conditions relevant to the study of generic origami patterns (finitely) folded out of a periodic crease pattern, e.g., Huffman grids, Ron-Resch pattern and Yoshimura pattern~\cite{Tachi2015}. Here is that stronger version.

\begin{lemma}
    Let $T_1$ and $T_2$ be positive real numbers and let $\ten A_1$ and $\ten A_2$ be two unitary linear maps of $\R^3$. Let $\ten x$ be a \emph{quasi-periodic} surface, i.e., such that
    \[
        \ten x_\alpha(\xi_1+m T_1,\xi_2+nT_2) = 
        \ten A_1^m\ten A_2^n\ten x_\alpha(\xi_1,\xi_2).
    \]
    Then,
    \[
        \int\scalar{\gt\omega,\dar_{\ten x}\ten w} = \int\scalar{\ten w,\dar_{\ten x}\gt\omega},
    \]
    for any $\gt\omega$ and $\ten w$ that are admissible and quasi-periodic.
\end{lemma}
\begin{remark}
    The average remains well-defined since $\scalar{\ten M\ten a,(\ten M\ten b)\wedge(\ten M\ten c)} = \scalar{\ten a,\ten b\wedge\ten c}$ for any unitary linear map $\ten M$.
\end{remark}
\begin{proof}
    Same as Lemma~\ref{lem:self} along with the above remark.
\end{proof}

\end{document}